\documentclass[twoside, 12pt]{article}
\usepackage{amssymb, amsmath, mathrsfs, amsthm}
\usepackage{graphicx}
\usepackage{color}
\usepackage[top=2cm, bottom=2cm, left=2cm, right=2cm]{geometry}
\usepackage{float, caption, subcaption}
\usepackage{ifpdf}
\usepackage{epstopdf}
\input{epsf}

\newcommand{\bd}{\begin{description}}
\newcommand{\ed}{\end{description}}
\newcommand{\bi}{\begin{itemize}}
\newcommand{\ei}{\end{itemize}}
\newcommand{\be}{\begin{enumerate}}
\newcommand{\ee}{\end{enumerate}}
\newcommand{\beq}{\begin{equation}}
\newcommand{\eeq}{\end{equation}}
\newcommand{\beqs}{\begin{eqnarray*}}
\newcommand{\eeqs}{\end{eqnarray*}}

\definecolor{DarkGreen}{rgb}{0.2, 0.6, 0.3}

\newtheorem{theorem}{Theorem}[section]

\newtheorem{definition}{Definition}
\newtheorem{corollary}[theorem]{Corollary}
\newtheorem{case}{Case}

\newtheorem{remark}{Remark}[section]

\begin{document}
\title{\textbf{Complete bipartite graphs without small rainbow stars} \footnote{Supported by the National
Science Foundation of China (No. 12061059) and the Qinghai Key Laboratory of
Internet of Things Project (2017-ZJ-Y21).} }

\author{Weizhen Chen\footnote{School of Mathematics and Statistis, Qinghai Normal University, Xining, Qinghai 810008, China. {\tt chenweizhenijk@163.com}}, \ \
Meng Ji\footnote{Corresponding author: College of Mathematical Science, Tianjin Normal University, Tianjin, China {\tt
jimengecho@163.com}}, \ \
Yaping Mao\footnote{Faculty of Environment and Information
Sciences, Yokohama National University, 79-2 Tokiwadai, Hodogaya-ku,
Yokohama 240-8501, Japan. {\tt mao-yaping-ht@ynu.ac.jp}}, \ \
Meiqin Wei\footnote{College of Arts and Sciences, Shanghai Maritime University, Shanghai, 201306, China. {\tt weimeiqin8912@163.com}}}
\date{}
\maketitle

\begin{abstract}
The $k$-edge-colored bipartite Gallai-Ramsey number $\operatorname{bgr}_k(G:H)$
is defined as the minimum integer $n$ such that $n^2\geq k$ and
for every $N\geq n$, every edge-coloring (using all $k$ colors) of complete bipartite graph $K_{N,N}$ contains a rainbow copy of $G$ or a monochromatic copy of $H$. In this paper, we first study the structural theorem on the complete bipartite graph $K_{n,n}$ with no rainbow copy of $K_{1,3}$. Next, we utilize the results to prove the exact values of $\operatorname{bgr}_{k}(P_4: H)$, $\operatorname{bgr}_{k}(P_5: H)$, $\operatorname{bgr}_{k}(K_{1,3}: H)$, where $H$ is a various union of cycles and paths and stars. \\[2mm]
\end{abstract}

\section{Introduction}

All graphs in this paper are undirected, finite and simple, and any undefined concepts or notations can be found in \cite{BM}. Let $G = (V(G),E(G))$ be a graph with vertex set $V(G)$ and edge set $E(G)$. Let $P_n$ and $C_n$ denote the path and cycle on $n$ vertices, respectively. For two vertex-sets $S$ and $T$, we denote all the edges one end in $S$ and the other in $T$ by $E(S,T)$.

A $k$-edge-coloring is \emph{exact} if all colors are used at least
once. In this work, we consider only exact edge-colorings of graphs. An edge-coloring of a graph is called \emph{rainbow} if no two edges have the same color. An edge-coloring of a graph is called \emph{monochromatic} if all the edges have the same color. For a $k$-edge-coloring of $G$ and a vertex subset $V'\subset V(G)$, we say \emph{$V'$ are incident with $t$ colors} if all the edges incident with $V'$ are colored by exactly $t$ colors.
Let $E^{(i)}$ be the set of edges with color $i$. The color degree $d_{c}(v)$ of a vertex $v$ is the number of different colors on edges incident with $v$. The \emph{color degree} $d_{c}(V')$ of a vertex-set $V'\subset V(G)$ is the number of different colors on edges between $V'$ and $V(G)\setminus V'$.
Let $\delta_{c}(G)$ and $\Delta_{c}(G)$ be minimum color degree and maximum color degree of the edge-colored graph $G$, respectively. For an edge $uv$, let $c(uv)$ denote the color used on $uv$. For a graph $G$, we use $rG$ to denote the union of $r$ disjoint copies of $G$. For graphs $G_{1}$ and $G_{2}$, we denote the union of $G_{1}$ and $G_{2}$ by $G_{1}\cup G_{2}$. For a complete bipartite graph, we also use $G(U, V)$ to denote the complete bipartite graph with bipartition $(U, V)$, where $|U|=|V|=n$. We denote a $k$-partition of $U$ (resp.,$V$) by $(U_{1},U_{2},\ldots,U_{k})$ (resp., $(V_{1},V_{2},\ldots,V_{k})$). And $c(U, V_i)=i$ means that all edges between $U$ and $V_{i}$ are colored by $i$.

Edge-colorings of complete graphs that contain no rainbow triangle have very interesting and somewhat
surprising structure. In 1967, Gallai \cite{Gallai} first examined this structure under the guise of transitive orientations. The result was proved again in \cite{GyarfasSimonyi} in the terminology of graphs.
\begin{theorem}{\upshape \cite{CameronEdmonds,Gallai,GyarfasSimonyi}}\label{Thm:G-Part}
In any edge-coloring of a complete graph containing no rainbow triangle,
there exists a nontrivial partition of the vertices (called a Gallai
partition), say $H_1,H_2,\ldots,H_t$, satisfying the following two
conditions.

$(a)$ The number of colors on the edges among $H_1,H_2,\ldots,H_t$
are at most two.

$(b)$ For each part pair $H_i,H_j \ (1\leq i\neq j\leq t)$, all the
edges between $H_i$ and $H_j$ receive the same color.
\end{theorem}

In 2010, Faudree et al \cite{FaudreeGouldJacobsonMagnant} and Gy\'{a}rf\'{a}s et al. \cite{GyarfasSarkozySeboSelkow} introduced the concept of Gallai-Ramsey numbers.

\begin{definition}{\upshape \cite{FaudreeGouldJacobsonMagnant, GyarfasSarkozySeboSelkow}}
Given two graphs $G$ and $H$, the general $k$-colored \emph{Gallai-Ramsey number} $\operatorname{gr}_{k}(G:H)$
is the minimum integer $n$ such that every $k$-edge-coloring of complete graph on $n$ vertices
contains either a rainbow copy of $G$ or a monochromatic copy of $H$.
\end{definition}

More results about Gallai-edge-colorings and Gallai-Ramsey numbers can be found in \cite{FujitaMagnant1,FujitaMagnant2,GyarfasLehelSchelp,GLST87,GyarfasSarkozySeboSelkow,
Hall-Magnant-Ozeki-Tsugaki,LBW22, LW20,Thomason-Wagner,Wang-Li,WangMaoZouMaganant, ZW21, ZW20}. And a survey can be found in \cite{FujitaMagnantOzeki}.

Motivated by giving structural theorems like Theorem \ref{Thm:G-Part} in Ramsey theory,
Thomason and Wagner \cite{Thomason-Wagner} studied the edge-colorings of complete graph $K_n$ that contains no rainbow path $P_{t+1}$ of length $t$. Schlage-Puchta and Wagner \cite{SW19} investigated the structural theorem of edge-colorings of complete graphs with no rainbow tree. Fujita and Magnant \cite{FujitaMagnant2} described the structure of rainbow $S_3^+$-free edge-colorings of a
complete graph, where the graph $S_3^+$ consisting of a triangle with a pendant edge.

Li et al. \cite{LiWangLiu19} studied the structure of complete bipartite graphs without rainbow paths $P_4$ and $P_5$, and we will use these results to prove our main results.
\begin{theorem}{\upshape \cite{LiWangLiu19}}\label{th-P4-Structure}
Let $K_{n,n}=G(U,V)$, $n\geq2$, be edge-colored such that it contains no rainbow $P_4$. Then one of the following propositions holds (replacing $U$ and $V$ if necessary):
\begin{description}
  \item[(a)] at most two colors are used;
  \item[(b)] $U$ can be partitioned into $k$ non-empty $U_1,U_2,\ldots, U_k$
such that $c(U_i,V)=i$ for $i=1, 2, \ldots, k$, where $k$ is the number of
colors used in the edge-coloring.
\end{description}
\end{theorem}

\begin{theorem}{\upshape \cite{LiWangLiu19}}\label{th-P5-Structure}
Let $K_{n,n}=G(U,V)$, $n\geq3$, be edge-colored such that it contains no rainbow $P_5$. Then, after renumbering the color, one of the following propositions holds (replacing $U$ and $V$ if necessary):
\begin{description}
  \item[(a)] at most three colors are used;
  \item[(b)] $U$ can be partitioned into two parts $U_1$ and $U_2$ with $|U_1|\geq 1$, $|U_2|\geq0$, and $V$ can be partitioned into $k$ parts $V_1, V_2, \ldots, V_{k}$ with $|V_1|\geq 0$, $|V_{j}|\geq1$ for $j=2, 3, \ldots, k$ such that $c(V_i, U_1)=i$, $c(V_i, U_2)=1$ for $i=1,2,\ldots, k$, where $k$ is the number of colors used in the edge-coloring;

  \item[(c)] $U$ can be partitioned into $k$ parts $U_1, U_2, \ldots, U_k$ with $|U_1|\geq 0$, $|U_{j}|\geq 1$, $j=2, 3, \ldots, k$, and $V$ can be partitioned into $k$ parts $V_1, V_2, \ldots V_{k}$ with $|V_1|\geq 0$, $|V_{j}|\geq 1$, $j=2, 3, \ldots, k$, such that only colors $1$ and $i$ can be used on the edges of $E(U_i, V_i)$, $i=1, 2, \ldots, k$, and every other edge is in $E^{(1)}$, where $k$ is the number of colors used in the edge-coloring;

  \item[(d)] $n=3$, $U=\{u_1, u_2, u_3\}$, $V=\{v_1, v_2, v_3\}$, $E^{(1)}=\{u_1v_1, u_2v_3, u_3v_2\}$, $E^{(2)}=\{u_1v_3, u_2v_1\}$, $E^{(3)}=\{u_1v_2, u_3v_1\}$ and $E^{(4)}=\{u_2v_2, u_3v_3\}$;
  \item[(e)] $n=4$, $U=\{u_1, u_2, u_3, u_4\}$, $V=\{v_1, v_2, v_3, v_4\}$, $E^{(1)}=\{u_1v_1, u_2v_3, u_3v_2, u_4v_4\}$, $E^{(2)}=\{u_1v_3, u_2v_1, u_3v_4, u_4v_2\}$, $E^{(3)}=\{u_1v_2, u_2v_4, u_3v_1, u_4v_3\}$ and $E^{(4)}=\{u_1v_4, u_2v_2, u_3v_3, u_4v_1\}$.
\end{description}
\end{theorem}
Li and Wang \cite{Li-Wang} obtained the following result, and we will generalize the result in Section $2$.
\begin{theorem}{\upshape \cite{Li-Wang}}
Let $K_{n,n} = G(U, V)$ for $n\geq3$ be an edge-colored complete bipartite graph with
no rainbow $K_{1,3}$. Then, after renumbering the colors, one of the following
propositions holds (replacing U and V if necessary):
\begin{description}
  \item[(a)] at most four colors are used;
  \item[(b)] $U$ can be partitioned into $k$ parts $U_{1}, U_{2},\ldots, U_{k}$ with $|U_{1}|\geq 0$, $|U_{j}|\geq1$, $j = 2, 3, \dots, k$, and $V$ can be partitioned into $k$ parts $V_{1}, V_{2}, \ldots , V_{k}$ with $|V_{1}|\geq 0$, $|V_{j}|\geq 1$, $j = 2, 3, \ldots , k$, such that only colors $1$ and $i$ can be used on the edges of $E(U_{i}, V_{i})$, $i = 1, 2, \ldots , k$, and every other edge is in
$E^{(1)}$, where $k$ is the number of colors used in the coloring.
\end{description}
\end{theorem}

Motivated by giving structures like Theorems \ref{th-P4-Structure} and \ref{th-P5-Structure}, we derive corresponding results for the $3$-star in Section $2$. In the following, we first introduce the concepts of bipartite Ramsey numbers defined by Faudree and Schelp \cite{FaudreeSchelp} and bipartite Gallai-Ramsey numbers defined by Li, Wang and Liu \cite{LiWangLiu19}.

\begin{definition}{\upshape \cite{FaudreeSchelp}}
The \emph{$k$-colored bipartite Ramsey number} $\operatorname{br}_k(H)$ is the minimum integer $n$ such that $n^2\geq k$ and for every $N\geq n$, every exact $k$-edge-coloring of the complete bipartite graph $K_{N,N}$ contains a monochromatic copy of $H$.
\end{definition}

\begin{definition}{\upshape \cite{LiWangLiu19}}
The \emph{$k$-colored bipartite Gallai-Ramsey number} $\operatorname{bgr}_k(G:H)$ is the minimum integer $n$ such that $n^2\geq k$ and for every $N\geq n$, every exact $k$-edge-coloring of complete bipartite graph $K_{N,N}$ contains either a rainbow copy of $G$ or a monochromatic copy of $H$.
\end{definition}

Clearly, for any bipartite graph $H$, we have that $\operatorname{bgr}_k(G:H)\leq\operatorname{br}_k(H)$. Some results on bipartite Ramsey number $\operatorname{br}_k(H)$ can be found in \cite{ErohOellermann,FaudreeSchelp,Hall-Magnant-Ozeki-Tsugaki,Wong}. Faudree-Schelp \cite{FaudreeSchelp} and Gy\'{a}rf\'{a}s-Lehel \cite{Gyarfas-Lehel} obtained the bipartite Ramsey numbers for $H=P_{n}$, respectively, which will be used in the proofs of our main results.
\begin{theorem}{\upshape \cite{FaudreeSchelp, Gyarfas-Lehel}}\label{br2-P_n}
For $n\geq 3$, we have
$$
\operatorname{br}_{2}(P_n)=\left\{
  \begin{array}{ll}
    n-1, & \hbox{$n$ \ is \ even;} \\
    n, & \hbox{$n$ \ is \ odd.}
  \end{array}
\right.$$
\end{theorem}

Li gave the upper and lower bounds for $\operatorname{bgr}_k(P_{4}:H)$ in \cite{Xihe-Li} , and in this paper we will close the gap for several graphs. Some notations will be needed.
Let $s(H)=\min\{|S|:(S,T)$ be a bipartition of $H$, $|S|\leq |T|\}$, $t(H)=\max\{|T|:(S,T)$ be a bipartition of $H$, $|S|\leq |T|\}$, $s^{*}(H)=\max\{|S|:(S,T)$ be a bipartition of $H$, $|S|\leq |T|\}$, and $t^{*}(H)=\min\{|T|:(S,T)$ be a bipartition of $H$, $|S|\leq |T|\}$.

\begin{theorem}{\upshape \cite{Xihe-Li}}
For a bipartite graph $H$, if $k\geq 3$ and $s(H)\geq 2$, then
\begin{align*}
\max\{t^{*}(H),(s(H)-1)k+1\}&\leq \operatorname{bgr}_k(P_{4}:H)\\[0.2cm]
&\leq \min\{\max\{t(H),(s(H)-1)k+1\},\max\{t^{*}(H),(s^{*}(H)-1)k+1\}\}.
\end{align*}
In particular, if $k\geq 3$ and $s(H)=s^{*}(H)\geq 2$, then $\operatorname{bgr}_k(P_{4}:H)= \max\{t(H),(s(H)-1)k+1\}.$
\end{theorem}

In Section $3$, we obtain several exact values of $\operatorname{bgr}_{k}(P_4: H)$ and $\operatorname{bgr}_{k}(P_5: H)$, respectively, where $H$ is a various union of cycles and paths. In Section $4$, $\operatorname{bgr}_{k}(K_{1,3}: H)$ are investigated, where $H$ is a star, or a complete bipartite graph.

\section{Proof of the structural theorem}

We give the structural theorem for rainbow $3$-star.
\begin{theorem}\label{th-K1,3-Structure}
Let $K_{n,n}=G(U,V)$ for $n\geq3$ be an edge colored complete bipartite graph with
no rainbow $K_{1,3}$. Then, after renumbering the colors, one of the
following propositions holds (replacing $U$ and $V$ if necessary):
\begin{description}
  \item[(a)] at most two colors are used;
  \item[(b)] $U$ can be partitioned into $k$ parts $U_1,U_2,\ldots,U_k$ with $|U_1|\geq0$, $|U_{i}|\geq1$, $i=2, \ldots k$, and $V$ can be partitioned into $k$ parts $V_1, V_2, \ldots, V_{k}$ with $|V_{1}|\geq 0$, $|V_{i}|\geq 1$, $i=2, \ldots k$, such that only colors $1$ and $i$ can be used on the edges of $E(U_i, V_i)$, $i=1, 2, \ldots k$, and every other edge is colored by color $1$, where $k$ is the number of colors used in the edge-coloring;
  \item[(c)] $U$ can be partitioned into $3$ parts $U_{1},U_{2},U_{3}$ such that each part is incident with at most two colors from $\{1,2,3\}$, and $V$ can be partitioned into $3$ parts $V_{1}, V_{2}, V_{3}$ such that each part is incident with at most two colors from $\{1,2,3\}$;
  \item[(d)] $U$ can be partitioned into $3$ parts $U_{1},U_{2},U_{3}$ such that each part is incident with two colors from $\{1,2\},\{2,3\}$ or $\{1,4\}$, and $V$ can be partitioned into $3$ parts $V_{1}, V_{2}, V_{3}$ such that each part is incident with two colors from $\{1,2\},\{1,3\}$ or $\{2,4\}$.
  \item[(e)] $U$ can be partitioned into $2$ parts $U_{1}$ and $U_{2}$, and $V$ can be partitioned into $4$ parts $V_{1}, V_{2}, V_{3}$ and $V_{4}$ (it is possible that $V_{i}=\emptyset$ for some $i\in [4]$) such that $c(U_{1},V_{1}\cup V_{2})=1$, $c(U_{1},V_{3}\cup V_{4})=4$, $c(U_{2}, V_{1}\cup V_{3})=2$ and $c(U_{2}, V_{2}\cup V_{4})=3$.
\end{description}
\end{theorem}
\begin{proof}
Let $c$ be an arbitrary exact $k$-edge-coloring of $K_{n,n}=G(U,V)$ without a rainbow copy of $K_{1,3}$. Then we have $\delta_{c}(G)\leq \Delta_{c}(G)\leq 2$. If $k\leq 2$, then $(a)$ holds. If $k=3$, then $(c)$ holds ($(b)$ is a special case of $(c)$ for $k=3$). In the following, we assume that $k\geq 4$.

For convenience, denoted by $C(v)$ the set of all colors incident with vertex $v$. If $\delta_{c}(G)=1$, say $C(u)=\{1\}$ for some $u\in U$, then we have $|C(v)\setminus\{1\}|\leq 1$ for any $v\in V$. If $|C(u')\setminus\{1\}|\leq 1$ for all $u'\in U\setminus \{u\}$, then $(b)$ holds; if $|C(v)\setminus\{1\}|\geq 2$ for some $u'\in U\setminus \{u\}$ (so $1\notin C(u')$ and $|C(u')|=2$), then $k\leq 3$ for avoiding a rainbow copy of $K_{1,3}$ centered at some $v\in V$, contradicting our assumption $k\geq 4$. Now we discuss the following two cases with the assumption that $\delta_{c}(G)= \Delta_{c}(G)= 2$.
\setcounter{case}{0}
\begin{case}
There exist two vertices $u,u'\in U$ with $C(u)\cap C(u')=\emptyset$.
\end{case}
Let $C(u_{1})=\{1,4\}$ and $C(u_{2})=\{2,3\}$. Then we have $C(u)\subset\{1,2,3,4\}$ for all $u\in U\setminus\{u_{1},u_{2}\}$. Otherwise, a rainbow copy of $K_{1,3}$ centered at some $v\in V$ will occur. If $C(u)\in \{\{1,4\},\{2,3\}\}$ for each $u\in U\setminus\{u_{1},u_{2}\}$, then easy to check that $(e)$ holds. If $C(u)\notin \{\{1,4\},\{2,3\}\}$ for some $u\in U\setminus\{u_{1},u_{2}\}$, say $C(u_{3})=\{1,2\}$, then $C(v)\in\{\{1,2\},\{1,3\},\{2,4\}\}$ with $k=4$ for each $v\in V$, otherwise, there is a rainbow copy of $K_{1,3}$ centered at some $v\in V$. Moreover, there must be some vertices $v_{1},v_{2}\in V$ with $C(v_{1})=\{1,3\}$ and $C(v_{2})=\{2,4\}$. If there is no vertex $v\in V$ with $C(v)=\{1,2\}$, then $(e)$ holds after replacing $U$ and $V$. If there is a vertex $v_{3}\in V $ with $C(v_{3})=\{1,2\}$, then $C(u)\in \{\{1,2\},\{1,4\},\{2,3\}\}$ for each $u\in U$, and so $(d)$ holds.

\begin{case}
$C(u)\cap C(u')\neq\emptyset$ for any two vertices $u,u'\in U$.
\end{case}

Without loss of generality, let $C(u_{1})=\{1,2\}$ and $C(u_{2})=\{1,3\}$ for $u_{1},u_{2}\in U$. Let $i$ be an arbitrary color in $[k]\setminus \{1,2,3\}$ and $u_{l}\in U$ be a vertex with $i\in C(u_{l})$. Then $l\geq 3$. Since $C(u_{1})\cap C(u_{l})\neq \emptyset$ and $C(u_{2})\cap C(u_{l})\neq \emptyset$, we have that $C(u_{l})=\{1,i\}$. Then $|C(v)\setminus \{1\}|\leq 1$ for each $v\in V$. It is easy to check that $(b)$ holds.

\end{proof}

\section{Results involving a rainbow path}\label{section-rainbow-P4P5}

In this section, we will show several exact values of bipartite Gallai-Ramsey numbers $\operatorname{bgr}_k(G:H)$ using the above characterizations.
\begin{theorem}
For integer $k\geq 3$, we have
$$
\operatorname{bgr}_{k}(P_4: H)=\left\{
  \begin{array}{ll}
    k(\lfloor r/2\rfloor+\lfloor(r+l)/2\rfloor)-k+1, & \hbox{$H= P_r\cup
P_{r+l}$\ for\ $r,l\geq2$;} \\[0.2cm]
    k(l_1+l_2)-k+1, & \hbox{$H= C_{2l_1}\cup C_{2l_2}$\ for\ $l_1, l_2\geq2$;}  \\[0.2cm]
    k(l+\lfloor m/2\rfloor)-k+1, & \hbox{$H=P_m\cup C_{2l}$\ for\ $m,l\geq2$.}
  \end{array}
\right.$$
\end{theorem}

\begin{proof}
Suppose that $H= P_r\cup
P_{r+l}$ ($H= C_{2l_1}\cup C_{2l_2}$, or $H=P_m\cup C_{2l}$). For the lower bound, we consider an exact $k$-edge-coloring of $K_{n,n}=G(U,V)$, where $n=k(\lfloor r/2\rfloor+\lfloor(r+l)/2\rfloor-1)$ ($n=k(l_1+l_2-1)$ and $|U_i|=(l_1+l_2-1)$, or $n=k(l+\lfloor m/2\rfloor)-k$ and $|U_i|=l+\lfloor m/2\rfloor-1$). Let $U=(U_1, U_2,\ldots,U_k)$ with
$|U_i|=\lfloor r/2\rfloor+\lfloor(r+l)/2\rfloor-1$ for $1\leq i\leq k$. We color the edges of $G$ such that $c(U_i, V)=i$ for $1\leq i\leq k$. It is easy to check that $G$ has neither
a rainbow copy of $P_4$ nor a monochromatic copy of $P_r\cup P_{r+l}$ ($C_{2l_1}\cup C_{2l_2}$ or $P_{m}\cup C_{2l}$).

For the upper bound, let $c$ be any $k$-edge-coloring of $K_{n,n}$ without rainbow copy of $P_4$, and let $G(U,V)$ be the $k$-edge-colored $K_{n,n}$ under the coloring $c$, where $n\geq k(\lfloor r/2\rfloor+\lfloor (r+l)/2\rfloor)+1-k$ ( $n\geq k(l_1+l_2)+1-k$, or $n\geq k(l+\lfloor m/2\rfloor)+1-k$). By Theorem \ref{th-P4-Structure}, we have a partition of $U$, say $U_1, U_2, \ldots, U_k$, such that $c(U_i, V)=i$ for $1\leq i\leq k$. Since
$$
\max\{|U_1|, |U_2|, \ldots,|U_k|\}\geq\lceil (k(\lfloor r/2\rfloor+\lfloor (r+l)/2\rfloor)+1-k)/k\rceil=(\lfloor r/2\rfloor+\lfloor (r+l)/2\rfloor)
$$
($\max\{|U_1|, |U_2|, \ldots,|U_k|\}\geq\lceil[k(l_1+l_2)+1-k]/k\rceil=(l_1+l_2)$, or $\max\{|U_1|, |U_2|, \ldots,|U_k|\}\geq\lceil (k(l+\lfloor m/2\rfloor)+1-k)/k\rceil=(l+\lfloor m/2\rfloor)$), it follows that $G(U_i,V)$ contains a monochromatic copy of $P_r\cup P_{r+l}$ ($C_{2l_1}\cup C_{2l_2}$ or $P_{m}\cup C_{2l}$) for $1\leq i\leq k$.
\end{proof}

In the following, we prove the exact values of bipartite Gallai-Ramsey numbers involving a rainbow $P_5$.

\begin{theorem}\label{P_5-rP_l}
For positive integers $k\geq4,l\geq 10$ and $r\geq 1$, we have
$$
\operatorname{bgr}_k(P_5: rP_l)=kr\lfloor l/2\rfloor-k+1.
$$
\end{theorem}

\begin{proof}
For the lower bound, we
consider a $k$-edge-coloring $c'$ of $K_{t,t}=G(U,V)$ with $t=kr\lfloor l/2\rfloor-k$.
We partition $V$ into $V_1, \ldots, V_{k}$ with $|V_i|=r\lfloor l/2\rfloor-1$ for $1\leq i\leq k$. Let $c'(U, V_i)=i$ for $1\leq i\leq k$. It is easy to
check that $K_{n,n}$ contains neither a rainbow copy of $P_5$ nor a monochromatic copy of $rP_l$. Hence, $\operatorname{bgr}_k(P_5: P_l)\geq kr\lfloor l/2\rfloor+1-k.$

Let $n\geq kr\lfloor l/2\rfloor+1-k$, and let $c$ be any
$k$-edge-coloring of $K_{n,n}$ without rainbow copy of $P_5$, and let $G(U,V)$ be the $k$-edge-colored $K_{n,n}$ under the coloring $c$. Suppose that there is no monochromatic copy of $rP_{l}$. We distinguish the following two cases to prove this theorem.
\setcounter{case}{0}
\begin{case}
It admits Theorem \ref{th-P5-Structure} $(b)$ for the edge-coloring $c$.
\end{case}

If $0\leq|U_2|\leq r\lfloor l/2\rfloor$, then
$|U_1| \geq rl-r\lfloor l/2\rfloor\geq r(l-\lfloor l/2\rfloor)\geq r\lceil l/2\rceil$.
Since
$$\max\{|V_1|, |V_2|, \ldots, |V_{k}|\}\geq \lceil (kr\lfloor l/2\rfloor+1-k)/k\rceil = r\lfloor l/2\rfloor,$$
then there exists a
monochromatic copy of $rP_l$ colored by $i$.
If $|U_2|> r\lfloor l/2\rfloor$, then there is a monochromatic copy of $rP_l$ colored by $1$ in $G(U_2\cup V)$, a contradiction.

\begin{case}
It admits Theorem \ref{th-P5-Structure} $(c)$ for the edge-coloring $c$.
\end{case}

If $|U_1|\geq r\lfloor l/2\rfloor$, then there is a monochromatic copy of $rP_l$ colored by $1$ in $G(U_1, V)$, a contradiction. Hence, $|U_1|\leq r\lfloor l/2\rfloor-1$. Then there exists a part, say $U_{2}$, with $U_{2}=\max\{|U_i|\,|\, i=2,\ldots, k\}\geq r\lfloor l/2\rfloor$.
If $|V_2|\leq r\lfloor l/2\rfloor-1$,
then $|V\setminus V_2|\geq r\lceil l/2\rceil$.
So there exists a monochromatic copy of $rP_l$ colored by $1$ in $G(U_i, (V\setminus V_i))$, a contradiction. Thus, $|V_2|\geq r\lceil l/2\rceil.$
If $|U\setminus U_{2}|\geq r\lfloor l/2\rfloor$ or $|U\setminus V_{2}|\geq r\lfloor l/2\rfloor$, then there is a monochromatic copy of $rP_{l}$ with color $1$, a contradiction. It forces $|U\setminus U_{2}|\leq r\lfloor l/2\rfloor-1$ and $|U\setminus V_{2}|\leq r\lfloor l/2\rfloor-1$. And then, we have $|U_{2}|\geq kr\lfloor l/2\rfloor-k+1-(r\lceil l/2\rceil-1)\geq 2r\lceil l/2\rceil$ and $|V_{2}|\geq kr\lfloor l/2\rfloor-k+1-(r\lceil l/2\rceil-1)\geq 2r\lceil l/2\rceil$. By Theorem \ref{br2-P_n}, there is a monochromatic copy of $P_{rl}$ (and thus a copy of $rP_{l}$), a contradiction.
\end{proof}

Using similar analysis as Theorem \ref{P_5-rP_l} above, we have an extend result of Theorem \ref{P_5-rP_l}, and we omit the proof details here.

\begin{theorem}\label{P_5-P_l-P_m}
For three positive integers $k, l_{i}, r$ with $k\geq 5$, $l_{i}\geq10$, we have
$$
\operatorname{bgr}_k\left(P_5: \bigcup^{r}_{i=1} P_{l_{i}}\right)=k\left(\sum^{r}_{i=1}\lfloor l_{i}/2\rfloor\right)-k+1.
$$
\end{theorem}

Now we focus on the monochromatic cycles and the union for cycles, paths, and complete bipartite graphs.

\begin{theorem}\label{P5-PlCm}
Let $r, s$ be non-negatives which are not both equal to $0$.  For positive integers $l_{i}\geq 2, m_{j}\geq 2$ and $k\geq 2+\sum^{r}_{i=1}(\lfloor l_{i}/2\rfloor+1)+\sum^{s}_{j=1}m_{j}$, we have
$$
\operatorname{bgr}_k\left(P_5: \bigcup^{r}_{i=1} P_{l_{i}}\cup \bigcup^{s}_{j=1}C_{2m_{j}}\right)=k\left(\sum^{r}_{i=1}\lfloor l_{i}/2\rfloor+\sum^{s}_{j=1}m_{j}\right)-k+1.
$$
\end{theorem}

\begin{proof}
For the lower bound, let $n=k(\sum^{r}_{i=1}\lfloor l_{i}/2\rfloor+\sum^{s}_{j=1}m_{j})-k$ and $K_{n,n}=G(U,V)$. We partition $V$ into $k$ parts,
say $V_1, V_2,\ldots, V_{k}$, where $|V_i|=\sum^{r}_{i=1}\lfloor l_{i}/2\rfloor+\sum^{s}_{j=1}m_{j}-1$.
Consider a $k$-edge-coloring $c'$ of $K_{n,n}=G(U,V)$
such that $c'(U, V_i)=i$ for $1\leq i\leq k$.
It is easy to check that $K_{n,n}$ contains neither a rainbow copy of $P_5$
nor a monochromatic copy of $\bigcup^{r}_{i=1} P_{l_{i}}\cup \bigcup^{s}_{j=1}C_{2m_{j}}$ for this edge-coloring.

For the upper bound, let $n\geq k\left(\sum^{r}_{i=1}\lfloor l_{i}/2\rfloor+\sum^{s}_{j=1}m_{j}\right)-k+1$ and $c$ be an arbitrary rainbow $P_5$-free
$k$-edge-coloring of $K_{n,n}=G(U,V)$. Assume that there exists a $k$-edge-coloring of $K_{n,n}$ containg no monochromatic $\bigcup^{r}_{i=1} P_{l_{i}}\cup \bigcup^{s}_{j=1}C_{2m_{j}}$. We proceed the remaining proof by discussing the following two cases.
\setcounter{case}{0}
\begin{case}
It admits Theorem \ref{th-P5-Structure} $(b)$ for the edge-coloring $c$.
\end{case}

If $0\leq|U_2|< \sum^{r}_{i=1}\lfloor l_{i}/2\rfloor+\sum^{s}_{j=1}m_{j}$, then $|U_1| \geq \sum^{r}_{i=1}\lceil l_{i}/2\rceil+\sum^{s}_{j=1}m_{j}.$
Since
$$\max\{|V_1|, |V_2|, \ldots, |V_{k}|\}=|V'| \geq \sum^{r}_{i=1}\lfloor l_{i}/2\rfloor+\sum^{s}_{j=1}m_{j},$$
it follows that there exists a
monochromatic copy of $\bigcup^{r}_{i=1} P_{l_{i}}\cup \bigcup^{s}_{j=1}C_{2m_{j}}$ colored by $i$ in the complete bipartite subgraph $G(U_1, V)$.
Suppose that $|U_2|\geq \sum^{r}_{i=1}\lfloor l_{i}/2\rfloor+\sum^{s}_{j=1}m_{j}$. Then there is a monochromatic copy of $\bigcup^{r}_{i=1} P_{l_{i}}\cup \bigcup^{s}_{j=1}C_{2m_{j}}$ colored by $1$ in $G(U_2, V)$, a contradiction.
\begin{case}
It admits Theorem \ref{th-P5-Structure} $(c)$ for the edge-coloring $c$.
\end{case}

If $\max\{|U_i|\}=|U_1|\geq \sum^{r}_{i=1}\lfloor l_{i}/2\rfloor+\sum^{s}_{j=1}m_{j}$, then there is a monochromatic copy of $\bigcup^{r}_{i=1} P_{l_{i}}\cup \bigcup^{s}_{j=1}C_{2m_{j}}$
with color $1$ in $G(U_1, V)$. So we have $\max\{|U_i|\,|\, i=2,\ldots, k\}|\geq \sum^{r}_{i=1}\lfloor l_{i}/2\rfloor+\sum^{s}_{j=1}m_{j}$. Without loss of generality, let $|U_{i}|\geq \sum^{r}_{i=1}\lfloor l_{i}/2\rfloor+\sum^{s}_{j=1}m_{j}$.
Suppose that there exist edges colored by $1$ in $G(U_i, V_i)$.
If $|V_i|\leq \sum^{r}_{i=1}\lfloor l_{i}/2\rfloor+\sum^{s}_{j=1}m_{j}-1$ ($|V\setminus V_i|\geq\sum^{r}_{i=1}\lceil l_{i}/2\rceil+\sum^{s}_{j=1}m_{j}$),
then there is a monochromatic copy of $\bigcup^{r}_{i=1} P_{l_{i}}\cup \bigcup^{s}_{j=1}C_{2m_{j}}$ colored by $1$ in $G(U_i, (V\setminus V_i))$, a contradiction. Hence, we have $|V_i|\geq \sum^{r}_{i=1}\lfloor l_{i}/2\rfloor+\sum^{s}_{j=1}m_{j}-1$.

Since $k\geq 2+\sum^{r}_{i=1}(\lfloor l_{i}/2\rfloor+1)+\sum^{s}_{j=1}m_{j}$ and $|U_{j}|\geq 1$ for
$j=2, \ldots, k$, it follows that $|\sum^{k}_{j=2,j\neq i}U_{j}|\geq k-2\geq \sum^{r}_{i=1}(\lfloor l_{i}/2\rfloor+1)+\sum^{s}_{j=1}m_{j}$. So there is a monochromatic complete bipartite subgraph $G(\sum^{k}_{j=2,j\neq i}U_{j}, V_{i})$ which contains a subgraph $\bigcup^{r}_{i=1} P_{l_{i}}\cup \bigcup^{s}_{j=1}C_{2m_{j}}$, a contradiction. Complete the proof.
\end{proof}

Clearly, if $r=0$ for Theorem \ref{P5-PlCm}, then we have the following corollary.
\begin{corollary}\label{cor-plcm}
For positive integers $k$, $l_{i}\geq2$ with $k\geq 2+\sum_{i=1}^{r} l_{i}$, we have
$$
\operatorname{bgr}_k\left(P_5: \bigcup_{i=1}^{r}C_{2l_{i}}\right)=k\sum_{i=1}^{r} l_{i}-k+1.
$$
\end{corollary}

\begin{theorem}\label{P_5-K_{t,t}}
For positive integers $k\geq t+2$ and $t\geq 3$, we have
$$
\operatorname{bgr}_k(P_5: K_{t,t})=kt-k+1.
$$
\end{theorem}

\begin{proof}
Let $n=kt-k$.
Consider a $k$-edge-coloring $c'$ of $K_{n,n}=G(U,V)$ such that $c'(U, V_i)=i$ for $1\leq i\leq k$, and we partition $V$ into $V_1,\cdots, V_{k}$ with $|V_i|=t-1 \ (1\leq i\leq k)$. It is easy to
check that $K_{n,n}$ contains neither a rainbow copy of $P_5$ nor a monochromatic copy of $K_{t,t}$. Hence, $\operatorname{bgr}_k(P_5: K_{t,t})\geq kt+1-k.$

Let $n\geq kt+1-k$ and $c$ be an arbitrary rainbow $P_5$-free
$k$-edge-coloring of $K_{n,n}=G(U,V)$. Assume that there exists a $k$-edge-coloring of $K_{n,n}$ containing no monochromatic copy of $K_{t,t}$. We proceed the proof by discussing the following two cases.
\setcounter{case}{0}
\begin{case}
It admits Theorem \ref{th-P5-Structure} $(b)$ for the edge-coloring $c$.
\end{case}

If $0\leq|U_2|\leq t-1$, then
$|U_1|\geq (k-1)t+1-k>t$. So
$\max\{|V_1|, |V_2|, \ldots, |V_{k}|\}=|V'| \geq t$.
As a result, there exists a monochromatic copy of $K_{t,t}$ colored by $i$ in the complete bipartite subgraph $G(U_1, V)$. If $|U_2|\geq t$, then there is a monochromatic copy of $K_{t,t}$ colored by $1$ in $G(U_2, V)$, a contradiction.

\begin{case}
It admits Theorem \ref{th-P5-Structure} $(c)$ for the edge-coloring $c$.
\end{case}

Let $|U_{i}|=\max\{|U_{j}|:1\leq j\leq k\}.$ Then $|U_{i}|\geq \lceil(kt-k+1)/k\rceil=t$. Since $|V_{j}|\geq 1$ for each $2\leq j\leq k$, it follows that $\sum^{k}_{j=1,j\neq i}|V_{j}|\geq k-2\geq t$, and hence there is a monochromatic copy of $K_{t,t}$ in $G(U_{i},V\setminus V_{i})$, a contradiction.
\end{proof}

\begin{remark}
Since $\bigcup_{i=1}^{r}C_{2l_{i}}$ is a subgraph of $K_{t,t}$ with $t=\sum_{i=1}^{r}l_{i}$, Corollary \ref{cor-plcm} is also a corollary of Theorem \ref{P_5-K_{t,t}}.
\end{remark}

\section{Results involving a rainbow star}

In this section, we get the exact values of bipartite Gallai-Ramsey numbers $\operatorname{bgr}_{k}(K_{1,3}: K_{1,t})$ and $\operatorname{bgr}_{k}(K_{1,3}: K_{s,t})$ using Theorem \ref{th-K1,3-Structure}, respectively.
\begin{theorem}\label{k13-k1t}
For integers $t\geq5, t\leq k\leq t+1$, we have
$$
\operatorname{bgr}_{k}(K_{1,3}: K_{1,t})=t+1.
$$
\end{theorem}

\begin{proof}
For the lower bound, we consider the following $k$-edge-coloring of $K_{t,t}=G(U,V)$. We partition $U$ (resp., $V$) into $k-1$ non-empty parts $U_{1},\ldots,U_{k-1}$ (resp., $V_{1},\ldots,V_{k-1}$). Let $c(U_{i},V_{i})=i$ for all $i\in [k-1]$, and all the other edges are of color $k$. Then $G$ contains neither a rainbow copy of $K_{1,3}$ nor a monochromatic copy of $K_{1,t}$.

For the upper bound, let $c$ be a $k$-edge-coloring of $K_{n,n}=G(U,V)$, where $n\geq t+1$. Suppose that $G$ contains neither a rainbow copy of $K_{1,3}$ nor a monochromatic copy of $K_{1,t}$. Since $k\geq t \geq 5$ and by Theorem \ref{th-K1,3-Structure}, we have Theorem \ref{th-K1,3-Structure} $(b)$ holds. If $|U_{i}|=1$ for some $2\leq i\leq k$, then $|U\setminus U_{i}|\geq t$. So there is a monochromatic copy of $K_{1,t}$ in $G(U\setminus U_{i},V_{i})$, a contradiction. Hence, $|U_{i}|\geq 2$ for all $2\leq i\leq k$. Thus, $\sum_{i=3}^{k}|U_{i}|\geq 2(k-2)\geq 2(t-2)\geq t$. So there is a monochromatic copy of $K_{1,t}$ in $G(\cup_{i=3}^{k}U_{i}, V_{2})$, a contradiction.
\end{proof}

\begin{theorem}
For integers $2\leq s<t, t\geq 5, t+s-\lfloor t/2\rfloor\leq k\leq t+s$, we have
$$
\operatorname{bgr}_{k}(K_{1,3}: K_{s,t})=t+s.
$$
\end{theorem}

\begin{proof}
For the lower bound, Let $K_{n,n}=G(U,V)$ with $n=s+t-1$. We partition $U$ $(resp., V)$ into $k-1$ non-empty parts $U_{1},\ldots, U_{k-1}$ $(resp., V_{1},\ldots, V_{k-1})$ with $|U_{i}|=|V_{i}|\leq t-1$ for each $i\in[k-1]$. We color the edges such that $c(U_{i},V_{i})=i$ for each $i\in[k-1]$, and all the other edges are of color $k$. Then $G$ has neither a rainbow copy of $K_{1,3}$ nor a monochromatic copy of $K_{s,t}$ of color $i\in [k-1]$. If there is a monochromatic copy of $K_{s,t}$ of color $k$, say in $G(\bigcup_{i=1}^{m}U_{i},\bigcup_{j=m+1}^{k-1}V_{j})$ with $|\bigcup_{i=1}^{m}U_{i}|\geq s$ and $|\bigcup_{j=m+1}^{k-1}V_{j}|\geq t$, then we have $n\geq s+t$ since $|U_{i}|=|V_{i}|$ for each $i\in [k-1]$, a contradiction. Hence, there is no monochromatic $K_{s,t}$ of color $k$.

For the upper bound, let $G$ be a rainbow $K_{1,3}$-free edge-coloring (using all $k$ colors) graph of $K_{n,n}=G(U,V)$ with $n\geq t+s$. By Theorem \ref{th-K1,3-Structure} ($b$), we partition $U$ (resp.,$V$) into $U_{1},\ldots,U_{k}$ (resp.,$V_{1},\ldots,V_{k}$). Without loss of generality, let $|U_{1}|\leq |U_{2}|\leq\ldots\leq |U_{k}|$.
If $|U_1|\geq s$ or $|V_1|\geq s$, then there is a monochromatic copy of $K_{s,t}$ in $G(U_1, V)$ or $G(U, V_1)$. So we assume that $|U_1|\leq s-1$ and $|V_1|\leq s-1$. Let
$$\alpha=\sum_{j=2}^{s+1}|U_{j}|$$ and $$\beta=|U_{1}|+\sum_{j=s+2}^{k}|U_{j}|.$$
Clearly, we have that $\alpha\geq s$ since $|U_{j}|\geq 1$ for $2\leq j\leq k$. Let $k=s+t-i$ for $0\leq i\leq \lfloor t/2\rfloor$. If $|U_{1}|=0$, then $\beta=\sum_{j=s+2}^{k}|U_{j}|$, and hence
\begin{align*}
  \beta&=\sum_{j=s+2}^{k-i}|U_{j}|+\sum_{j=k-i+1}^{k}|U_{j}|\\[0.2cm]
&\geq \left(k-i-(s+2)+1\right)+\left((k-(k-i+1)+1)+i+1\right)\\[0.2cm]
&=k+i-s=t.
\end{align*}
If $|U_{1}|\geq 1$, then we have that
\begin{align*}
  \beta&=|U_{1}|+\sum_{j=s+2}^{k-i}|U_{j}|+\sum_{j=k-i+1}^{k}|U_{j}|\\[0.2cm]
&\geq 1+\left(k-i-(s+2)+1\right)+\left((k-(k-i+1)+1)+i\right)\\[0.2cm]
&=k+i-s=t,
\end{align*}
and hence $\alpha\geq s$ and $\beta\geq t$. If $|\sum_{j=2}^{s+1}V_{j}|\geq s$, then there is a monochromatic copy of $K_{s,t}$ in $G(U\setminus \bigcup_{j=2}^{s+1}U_{j},\sum_{j=2}^{s+1}V_{j})$; if $|\sum_{j=2}^{s+1}V_{j}|< s$, then there is a monochromatic copy of $K_{s,t}$ in $G( \bigcup_{j=2}^{s+1}U_{j},V\setminus\sum_{j=2}^{s+1}V_{j})$. Completed the proof.
\end{proof}

\vspace*{4mm}

\noindent {\bf Acknowledgment.} We sincerely thank both referees for their helpful comments to improve our presentation and simplify the proofs.
\end{document}